\documentclass[14pt, draft]{article}
\usepackage[cp1251]{inputenc}
\usepackage[ukrainian, russian, english]{babel}

\usepackage{ amsfonts, amssymb}
\usepackage{amsmath}
\usepackage{amsthm}

\begin{document}

\selectlanguage{ukrainian} \thispagestyle{empty}
 \pagestyle{myheadings}              

\pagestyle{myheadings}              

УДК 517.5 \vskip 3mm

\noindent \bf А.С. Сердюк  \rm (Інститут математики НАН України, Київ) \\
\noindent \bf Т.А. Степанюк  \rm (Східноєвропейський нац.
університет імені Лесі Українки, Луцьк)

\noindent {\bf A.S. Serdyuk} (Institute of Mathematics NAS of
Ukraine, Kyiv) \\
 \noindent {\bf T.A. Stepaniuk} (Lesya Ukrainka Eastern European National
 University, Lutsk)\\

 \vskip 5mm

 {\bf Наближення класів узагальнених інтегралів Пуассона сумами Фур'є в метриках просторів $L_s$}

\vskip 5mm

{\bf Approximations of  classes of generalized Poisson integrals by Fourier sums in metrics of spaces $L_{s}$}

\vskip 5mm

 \rm В метриках пространств $L_{s}, \ 1\leq s\leq\infty$, найдены асимптотические равенства для  верхних граней приближений суммами Фурье  на классах
  обобщенных интегралов Пуассона периодических функций, принадлежащих единичному шару  пространства $L_{1}$.

\vskip 5mm

 \rm In metrics of spaces $L_{s}, \ 1\leq s\leq\infty$, we find asymptotic equalities for upper bounds of approximations by Fourier  sums  on classes of generalized Poisson integrals of
  periodic functions,   which belong to unit ball of space $L_{1}$.
\newpage

Дана робота тісно пов'язана з роботами авторів [\ref{SerdyukStepanyuk2016dop}, \ref{SerdyukStepanyuk2016}]. В ній продовжуються дослідження апроксимативних властивостей сум Фур'є на класах узагальнених інтегралів Пуассона $C^{\alpha,r}_{\beta,p}$.

Нехай $L_{s}$,
$1\leq s<\infty$, --- простір $2\pi$--періодичних сумовних в $s$--ій
степені на $[0,2\pi)$ функцій $f$, в якому  норма задана формулою
${\|f\|_{s}=\Big(\int\limits_{0}^{2\pi}|f(t)|^{s}dt\Big)^{\frac{1}{s}}}$; $L_{\infty}$ --- простір
$2\pi$--періодичних вимірних і суттєво обмежених функцій $f$ з нормою
$\|f\|_{\infty}=\mathop{\rm{ess}\sup}\limits_{t}|f(t)|$; $C$ --- простір $2\pi$--періодичних неперервних функцій $f$, у
якому норма задається за допомогою рівності $
{\|f\|_{C}=\max\limits_{t}|f(t)|}$.

Позначимо через $C^{\alpha,r}_{\beta,p}, \ \alpha>0, \ r>0, \ 1\leq p\leq\infty,$ множину всіх $2\pi$--періодичних функцій, котрі  при всіх
$x\in\mathbb{R}$ зображуються за допомогою згортки (див., наприклад, [\ref{Stepanets1}, с.
144])
\begin{equation}\label{conv}
f(x)=\frac{a_{0}}{2}+\frac{1}{\pi}\int\limits_{-\pi}^{\pi}P_{\alpha,r,\beta}(x-t)\varphi(t)dt,
\ a_{0}\in\mathbb{R}, \ \varphi\in B_{p}^{0}, \
\end{equation}
$$
B_{p}^{0}=\left\{\varphi: \ ||\varphi||_{p}\leq 1, \  \varphi\perp1\right\},
 \ 1\leq p\leq \infty,
$$
з узагальненими ядрами Пуассона вигляду
\begin{equation}\label{kernel}
P_{\alpha,r,\beta}(t)=\sum\limits_{k=1}^{\infty}e^{-\alpha k^{r}}\cos
\big(kt-\frac{\beta\pi}{2}\big), \ \ \alpha>0, \ \ r>0, \ \ \beta\in
    \mathbb{R}.
\end{equation}
Функції $f$ вигляду (\ref{conv}) називають узагальненими інтегралами Пуассона функцій $\varphi$.  При $r=1$ класи $C^{\alpha,r}_{\beta,p}$ є відомими класами інтегралів Пуассона.

При довільних $r>0$ класи  $C^{\alpha,r}_{\beta,p}$ належать до множини $D^{\infty}$
 нескінченно диференційовних $2\pi$--періодичних функцій, тобто $C^{\alpha,r}_{\beta,p}\subset D^{\infty}$ (див., наприклад, [\ref{Stepanets1}, с. 139; \ref{Stepanets_Serdyuk_Shydlich}, с. 1408]). Більш того,  як випливає з теореми 1 роботи [\ref{Stepanets_Serdyuk_Shydlich2009}], при довільних $r>0$, $\alpha>0$,
 $\beta\in\mathbb{R}$ i $1\leq p\leq\infty$ має місце вкладення $C^{\alpha,r}_{\beta,p}\subset \mathcal{J}_{1/r}$, де $\mathcal{J}_{a}, a>0,$ --- відомі класи Жевре
$$
\mathcal{J}_{a}=\bigg\{f\in D^{\infty}: \ \sup\limits_{k\in \mathbb{N}}\Big(\frac{\|f^{(k)}\|_{C}}{(k!)^{a}}\Big)^{1/k}<\infty \bigg\}.
$$
 При $r\geq1$ класи $C^{\alpha,r}_{\beta,p}$
складаються з  функцій $f$, які  допускають регулярне продовження в смугу $|\mathrm{Im} \ z|\leq c, \ c>0$ комплексної площини (див., наприклад, [\ref{Stepanets1}, с.~141]), тобто є класами аналітичних функцій.
При $r>1$
класи  $C^{\alpha,r}_{\beta,p}$  складаються з функцій регулярних в усій комплексній площині, тобто є класами цілих функцій (див., наприклад, [\ref{Stepanets1}, с. 142]).

Розглянемо апроксимативні характеристики вигляду
    \begin{equation}\label{sum}
 {\cal E}_{n}(C^{\alpha,r}_{\beta,p})_{s}=\sup\limits_{f\in
C^{\alpha,r}_{\beta,p}}\|f(\cdot)-S_{n-1}(f;\cdot)\|_{s},  \ r>0, \ \alpha>0, \ \beta\in\mathbb{R}, \ 1\leq p,s \leq \infty,
  \end{equation}
  де $S_{n-1}(f;\cdot)$ --- частинні суми Фур'є порядку $n-1$ функції $f$.

  На даний час порядкові оцінки величин ${\cal E}_{n}(C^{\alpha,r}_{\beta,p})_{s}$ вигляду (\ref{sum}) відомі при всіх допустимих значеннях параметрів $\alpha$, $r$, $\beta$, $p$ i $s$ (див. роботи [\ref{Stepanets2}, с. 60; \ref{S_S}--\ref{SerdyukStepanyukBulleten}] і наявну в них бібліографію). У даній роботі досліджуватимуться величини (\ref{sum}) при $p=1$ i $r\in(0,1)$ з метою одержання для них асимптотичних при $n\rightarrow\infty$ рівностей.

  При $s=\infty$ асимптотичні рівності для величин вигляду (\ref{sum}) відомі при всіх $r>0$, $\alpha>0$, $\beta\in\mathbb{R}$ i $1\leq p \leq \infty$
  (див. роботи [\ref{SerdyukStepanyuk2016dop}, \ref{SerdyukStepanyuk2016}]). В зазначеному випадку норму $\|\cdot\|_\infty$ у (\ref{sum}) можна замінити на $\|\cdot\|_C$, тобто ${\cal E}_{n}(C^{\alpha,r}_{\beta,p})_{\infty}={\cal E}_{n}(C^{\alpha,r}_{\beta,p})_{C}$.

  Що ж стосується питання про сильну асимптотику при $n\rightarrow\infty$ величини (\ref{sum}) у випадку $p=1$ зазначимо наступне.

  Випадок ${r=1}$, $s=1$ дослідив С.М. Нікольський [\ref{Nikolsky 1946}, с. 221], встановивши асимптотичну рівність
 \begin{equation}\label{triang}
  {\cal E}_{n}(C^{\alpha,1}_{\beta,1})_{1}=e^{-\alpha n}\Big(\frac{8}{\pi^{2}}{\bf K}(e^{-\alpha})+O(1)n^{-1}\Big),
 \end{equation}
      де
  $$
\mathbf{K}(q):=\int\limits_{0}^{\frac{\pi}{2}}\frac{dt}{\sqrt{1-q^{2}\sin^{2}t}}, \ q\in(0,1),
$$
--- повний еліптичний інтеграл першого роду, а  $O(1)$ --- величина рівномірно обмежена відносно параметрів  $n$ i $\beta$.

  Згодом рівність (\ref{triang}) уточнив С.Б. Стєчкін в [\ref{Stechkin 1980}, с. 139],  встановивши асимптотичну формулу
\begin{equation}\label{stechkin}
  {\cal E}_{n}(C^{\alpha,1}_{\beta,1})_{1}=
 e^{-\alpha n}\Big(\frac{8}{\pi^{2}}\mathbf{K}(e^{-\alpha})+O(1)\frac{e^{-\alpha}}{(1-e^{-\alpha})n}\Big), \  \ \alpha>0, \ \beta\in\mathbb{R},
  \end{equation}
в якій величина  $O(1)$   рівномірно обмежена відносно  всіх розглядуваних параметрів.

В  роботі [\ref{Serdyuk2005Int}] при $r=1$ і довільних $1\leq s\leq\infty$ для величин ${\cal E}_{n}(C^{\alpha,r}_{\beta,1})_{s}$, $\alpha>0$, $\beta\in\mathbb{R}$,
встановлено  рівність
\begin{equation}\label{hhj}
  {\cal E}_{n}(C^{\alpha,1}_{\beta,1})_{s}=e^{-\alpha n}\bigg(\frac{2}{\pi^{1+\frac{1}{s}}}\|\cos t\|_{s}K(s,e^{-\alpha})+O(1)\frac{e^{-\alpha}}{n(1-e^{-\alpha})^{\sigma(s)}}\bigg), \ \
  \end{equation}
де
\begin{equation}\label{star1}
\sigma(s):={\left\{\begin{array}{cc}
1, \ & s=1,  \\
2, \ & s\in(1, \infty],
  \end{array} \right.}
 \end{equation}
\begin{equation}\label{star2}
K(s,q):=\frac{1}{2^{1+\frac{1}{s}}}\Big\|(1-2q\cos t+q^{2})^{-\frac{1}{2}} \Big\|_{s}, \ q\in(0,1),
 \end{equation}
a величина $O(1)$ рівномірно обмежена відносно $n$, $s$, $\alpha$ i $\beta$.
При $s=1$, з урахуванням  рівності $K(1,q)=\mathbf{K}(q)$, оцінка (\ref{hhj}) збігається з оцінкою  (\ref{stechkin}).

Як показано в [\ref{Serdyuk2012}], при $1\leq s< \infty$ має місце наступна рівність:
\begin{equation}\label{hyper2}
K(s,q)=\frac{\pi^{\frac{1}{s}}}{2}F^{\frac{1}{s}}\Big(\frac{s}{2}, \frac{s}{2}; 1; q^{2}\Big),  \ q\in(0,1),
\end{equation}
де $F(a,b;c;z)$ гіпергеометрична функція Гаусса
\begin{equation}\label{hyper}
F(a,b;c;z)=1+\sum\limits_{k=1}^{\infty}\frac{(a)_{k}(b)_{k}}{(c)_{k}}\frac{z^{k}}{k!},
\end{equation}
\begin{equation}\label{hyper1}
(x)_{k}:=x(x+1)(x+2)...(x+k-1).
\end{equation}

В [\ref{Serdyuk2013}, с. 250] було доведено, що для величин ${\cal E}_{n}(C^{\alpha,1}_{\beta,1})_{s}$ при $s=2$  виконується рівність
\begin{equation}\label{star3}
 {\cal E}_{n}(C^{\alpha,1}_{\beta,1})_{2}=\frac{1}{\sqrt{\pi(1-e^{-2\alpha})}}e^{-\alpha n}, \ \alpha>0, \ \beta\in\mathbb{R}, \ n\in\mathbb{N}.
\end{equation}

\noindent Рівність (\ref{star3}) уточнює асимптотичну рівність (\ref{hhj}) при $s=2$  в наступному сенсі: зазначена рівність (\ref{hhj}) при $s=2$ залишається вірною, якщо в ній обнулити залишковий член. Отже, взявши до уваги формули (\ref{star1}), (\ref{star2}) і (\ref{star3}) та очевидну рівність $F(1,1;1;q^{2})=\frac{1}{1-q^{2}}, \ \ q\in(0,1)$, для усіх $\alpha>0$, $\beta\in\mathbb{R}$, і $1\leq s\leq\infty$ можемо записати формулу
\begin{equation}\label{star4}
 {\cal E}_{n}(C^{\alpha,1}_{\beta,1})_{s}={\left\{\begin{array}{cc}
e^{-\alpha n}\bigg(\frac{\|\cos t\|_{s}}{\pi}F^{\frac{1}{s}}\Big(\frac{s}{2},\frac{s}{2};1;e^{-\alpha}\Big) +O(1)\frac{\xi(s)e^{-\alpha}}{n(1-e^{-\alpha})^{\sigma(s)}}\bigg), \ & 1\leq s<\infty,  \\
e^{-\alpha n}\bigg(\frac{1}{\pi(1-e^{-\alpha})}+O(1)\frac{e^{-\alpha}}{n(1-e^{-\alpha})^{2}}\bigg), \ \ \ \ \ \ \ \ \ \ \ \ \ \ \ \ \ \ \ \ \ \ & s=\infty,
  \end{array} \right.}
\end{equation}
в якій
$$
 \xi(s)={\left\{\begin{array}{cc}
0, \ & s=2,  \ \ \ \ \ \ \ \ \ \ \ \ \ \ \ \ \  \\
1, \ & s\in[1,2)\cup(2,\infty),
  \end{array} \right.}
$$
$\sigma(s)$ означена формулою (\ref{star1}), a $O(1)$ --- величина, що рівномірно обмежена відносно усіх розглядуваних параметрів.

Крім того, як випливає з [\ref{Serdyuk2013}], при  всіх $r>0$
 має місце рівність
\begin{equation}\label{serd2011}
 {\cal E}_{n}(C^{\alpha,r}_{\beta,1})_{2}=\frac{1}{\sqrt{\pi}}\Big(\sum\limits_{k=n}^{\infty}e^{-2\alpha k^{r}}\Big)^{\frac{1}{2}}, \ \alpha>0, \ r>0, \ \beta\in\mathbb{R}, \ n\in\mathbb{N}.
\end{equation}

 У випадку ${r>1}$ і  $s=1$ асимптотично точні оцінки для величин ${\cal E}_{n}(C^{\alpha,r}_{\beta,1})_{s}$, $\alpha>0$, $\beta\in\mathbb{R}$,   були одержані О.І. Степанцем
 [\ref{Step monog 1987}, с.155], який показав, що для довільних $n\in\mathbb{N}$
\begin{equation}\label{step_prepr}
  {\cal E}_{n}(C^{\alpha,r}_{\beta,1})_{1}=\Big(\frac{4}{\pi}+\gamma_{n}\Big)e^{-\alpha n^{r}},
\end{equation}
де
$$
|\gamma_{n}|<2\Big(1+\frac{1}{\alpha r n^{r-1}}\Big)e^{-\alpha r n^{r-1}}.
$$

 Згодом С.О. Теляковський [\ref{Teljakovsky 1989}, с. 517] при $r>1$, $\alpha>0$ i $\beta\in\mathbb{R}$ встановив  асимптотичну рівність
\begin{equation}\label{tel}
  {\cal E}_{n}(C^{\alpha,r}_{\beta,1})_{1}=\frac{4}{\pi}e^{-\alpha n^{r}}+
  O(1)\Big(e^{-\alpha ( 2(n+1)^{r}-n^{r})}
  +\Big(1+\frac{1}{\alpha r(n+2)^{r}}\Big)e^{-\alpha (n+2)^{r}}  \Big),
\end{equation}
де  $O(1)$ --- величина, що рівномірно обмежена відносно  всіх розглядуваних параметрів.
Формула (\ref{tel}) містить більш точну оцінку залишкового члена в асимптотичному розкладі величин ${\cal E}_{n}(C^{\alpha,r}_{\beta,1})_{1}$ у порівнянні з оцінкою (\ref{step_prepr}).

При $r>1$ і довільних $1\leq s\leq\infty$ асимптотичні рівності для величин ${\cal E}_{n}(C^{\alpha,r}_{\beta,1})_{s}$, $\alpha>0$, $\beta\in\mathbb{R}$,
знайдені в [\ref{Serdyuk2005Int}, с. 1408] і мають вигляд
\begin{equation}\label{ser}
  {\cal E}_{n}(C^{\alpha,r}_{\beta,1})_{s}=e^{-\alpha n^{r}}\Big(\frac{\|\cos t\|_{s}}{\pi}+O(1)\Big(1+\frac{1}{\alpha r n^{r-1}}\Big)e^{-\alpha rn^{r-1}}\Big),
\end{equation}
де величина $O(1)$ рівномірно обмежена відносно всіх розглядуваних параметрів. При $s=1$ формула (\ref{ser}) випливає із (\ref{step_prepr}) i
(\ref{tel}).

 Що ж стосується  випадку ${0<r<1}$,   то асимптотичні рівності для величин ${\cal E}_{n}(C^{\alpha,r}_{\beta,1})_{s}$,  $\alpha>0$,  $\beta\in\mathbb{R}$, за виключенням наведеного випадку $s=2$,  були відомі у випадку $s=1$ (див. [\ref{Step monog 1987}, с.153])
\begin{equation}\label{step}
  {\cal E}_{n}(C^{\alpha,r}_{\beta,1})_{1}=\frac{4}{\pi^{2}}e^{-\alpha n^{r}}\ln n^{1-r}+O(1)e^{-\alpha n^{r}},
\end{equation}
а також у випадку $s=\infty$ (див.  [\ref{SerdyukStepanyuk2016dop}, \ref{SerdyukStepanyuk2016}])
\begin{equation}\label{SerStep}
{\cal E}_{n}(C^{\alpha,r}_{\beta,1})_{\infty}=
e^{-\alpha n^{r}}n^{1-r}\bigg(\frac{1}{\pi\alpha r}+O(1)\Big(\frac{1}{(\alpha r)^{2}}\frac{1}{n^{r}}+\frac{1}{n^{1-r}}\Big)\bigg),
\end{equation}
де $O(1)$ --- величини, що рівномірно обмежені відносно  параметрів $n$ i $\beta$.

Зауважимо, що на класах $C^{\alpha,r}_{\beta,1}$,  $\alpha>0$,  $\beta\in\mathbb{R}$, при ${0<r<1}$ у випадку наближень в метриках просторів $L_{s}, {1< s\leq\infty}$,  суми Фур'є забезпечують порядок найкращих наближень тригонометричними поліномами, тобто (див., наприклад, [\ref{S_S}, \ref{S_S2}])
$$
{\cal E}_{n}(C^{\alpha,r}_{\beta,1})_{s}\asymp{ E}_{n}(C^{\alpha,r}_{\beta,1})_{s}\asymp e^{-\alpha n^{r}}n^{\frac{1-r}{s'}}, \ \ \frac{1}{s}+\frac{1}{s'}=1,
$$
де
$$
{ E}_{n}(C^{\alpha,r}_{\beta,1})_{s}=\sup\limits_{f\in
C^{\alpha,r}_{\beta,1}}\, \inf\limits_{t_{n-1}\in\mathcal{T}_{2n-1}}\|f-t_{n-1}\|_{s},
$$
 $\mathcal{T}_{2n-1}$ --- підпростір усіх тригонометричних
поліномів $t_{n-1}$ порядку не вищого за ${n-1}$.

В даній роботі встановлено асимптотично непокращувані оцінки величин ${\cal E}_{n}(C^{\alpha,r}_{\beta,1})_{s}$,  ${\alpha>0}$,  $\beta\in\mathbb{R}$,
 при довільних $0< r<1$ i $1\leq s\leq\infty$. При цьому  в отриманих асимптотичних рівностях в явному вигляді виписані
 оцінки залишкового члена через параметри задачі, що може бути корисним для практичного застосування.

При довільних фіксованих $\alpha>0$, $r\in(0,1)$ i $1\leq p\leq\infty$ через
$n_0=n_0(\alpha,r,p)$ позначимо найменший з номерів $n$ такий, що
\begin{equation}\label{n_p}
 \frac{1}{\alpha r}\frac{1}{n^{r}}+\frac{\alpha r \chi(p)}{n^{1-r}}\leq{\left\{\begin{array}{cc}
 \frac{1}{14},  & p=1, \\
\frac{1}{(3\pi)^3}\cdot\frac{p-1}{p}, & 1< p<\infty, \\
\frac{1}{(3\pi)^3}, & p=\infty, \
  \end{array} \right.}
\end{equation}
де $\chi(p)=p$ при $1 \leq p<\infty$ i $\chi(p)=1$ при $p=\infty$; через
 $n_1=n_1(\alpha,r)$ --- найменший з номерів $n$ такий, що
\begin{equation}\label{n_2}
 \frac{1}{\alpha r}\frac{1}{n^{r}}\Big(1+\ln\Big(\frac{\pi  n^{1-r}}{\alpha r}\Big)\Big)+\frac{\alpha r }{n^{1-r}}\leq
\frac{1}{(3\pi)^3}.
\end{equation}

Також при довільних $\upsilon>0$ i ${1\leq s\leq \infty}$ покладемо
\begin{equation}\label{norm_j}
  \mathcal{J}_{s}(\upsilon):=\big\|\frac{1}{\sqrt{t^{2}+1}} \big\|_{L_{s}[0,\upsilon]},
\end{equation}
де
$$
\|f \|_{L_{s}[a,b]}=
   {\left\{\begin{array}{cc}
\bigg(\int\limits_{a}^{b}|f(t)|^{s}dt
\bigg)^{\frac{1}{s}}, & 1\leq s<\infty, \\
\mathop{\rm{ess}\sup}\limits_{t\in[a,b]}|f(t)|, \ & s=\infty. \
  \end{array} \right.}
$$

В прийнятих позначеннях має місце наступне твердження.

{\bf Теорема 1.} {\it Нехай  $0<r<1$, $1\leq s\leq\infty$, $\alpha>0$ i $\beta\in\mathbb{R}$. Тоді при   $n\geq n_0(\alpha,r,s')$, $\frac{1}{s}+\frac{1}{s'}=1$, справедлива  оцінка
\begin{equation}\label{theorem1}
{\cal E}_{n}(C^{\alpha,r}_{\beta,1})_{s}=
e^{-\alpha n^{r}}n^{\frac{1-r}{s'}}\bigg(\frac{\|\cos t\|_{s}}{\pi^{1+\frac{1}{s}}(\alpha r)^{\frac{1}{s'}}} \mathcal{J}_{s}\Big(\frac{ \pi n^{1-r}}{\alpha r}\Big)+\gamma_{n,s}^{(1)}\Big(\frac{1}{(\alpha r)^{1+\frac{1}{s'}}} \mathcal{J}_{s}\Big(\frac{ \pi n^{1-r}}{\alpha r}\Big)\frac{1}{n^{r}}+\frac{1}{n^{\frac{1-r}{s'}}}\Big)\bigg),
\end{equation}
де для величини ${\gamma_{n,s}^{(1)}=\gamma_{n,s}^{(1)}(\alpha,r,\beta)}$ виконується нерівність  ${|\gamma_{n,s}^{(1)}|\leq(14\pi)^{2}}$}.

{\bf Доведення теореми 1.}
Згідно з  (\ref{conv}) і (\ref{sum})
\begin{equation}\label{f1}
{\cal E}_{n}(C^{\alpha,r}_{\beta,1})_{s}=
\frac{1}{\pi}\sup\limits_{\varphi\in B_{1}^{0}}\bigg\|\int\limits_{-\pi}^{\pi}P_{\alpha,r,\beta}^{(n)}(x-t)\varphi(t)dt\bigg\|_{s}, \ \
  \ 1\leq s\leq\infty,
\end{equation}
 де
 \begin{equation}\label{f2}
P_{\alpha,r,\beta}^{(n)}(t):=
\sum\limits_{k=n}^{\infty}e^{-\alpha k^{r}}\cos\Big(kt-\frac{\beta\pi}{2}\Big),  \ 0<r<1, \ \alpha>0, \ \beta\in\mathbb{R}.
\end{equation}

Далі скористаємось наступним твердженням з роботи [\ref{Serdyuk2005Int}, с. 1398].

{\bf Лема $1$.} {\it Нехай $K\in L_{s}$, $1\leq s\leq\infty$. Тоді для величини
$$
\mathcal{E}(K)_{s}:=\sup\limits_{\varphi\in B_{1}^{0}}\bigg\|\frac{1}{\pi}\int\limits_{-\pi}^{\pi}\varphi(\cdot-t)K(t)dt\bigg\|_{s}
$$
виконуються співвідношення}
\begin{equation}\label{star5}
\frac{1}{2\pi}\sup\limits_{h\in\mathbb{R}}\|K(\cdot)-K(\cdot+h) \|_{s}\leq \mathcal{E}(K)_{s}\leq\frac{1}{\pi}\|K \|_{s}.
\end{equation}

Покладаючи в умовах леми $1$ $K(t)=P_{\alpha,r,\beta}^{(n)}(t)$ і враховуючи рівність (\ref{f1}), з (\ref{star5}) отримуємо співвідношення
\begin{equation}\label{s12}
 \frac{1}{2\pi}\sup\limits_{h\in\mathbb{R}}\|P_{\alpha,r,\beta}^{(n)}(\cdot)-P_{\alpha,r,\beta}^{(n)}(\cdot+h) \|_{s}\leq \mathcal{E}(C^{\alpha,r}_{\beta,1})_{s}\leq\frac{1}{\pi}\|P_{\alpha,r,\beta}^{(n)} \|_{s}, \ \ 1\leq s\leq\infty.
\end{equation}

Згідно з  [\ref{SerdyukStepanyuk2016dop}] та  [\ref{SerdyukStepanyuk2016}] при довільних $r\in(0,1)$, $\alpha>0$, $\beta\in\mathbb{R}$, $1\leq s\leq\infty$, $\frac{1}{s}+\frac{1}{s'}=1$ i $n\in\mathbb{N}$ при $n\geq n_0(\alpha,r,s')$
мають місце співвідношення
$$
\frac{1}{\pi}\|P_{\alpha,r,\beta}^{(n)} \|_{s}= \ \ \ \ \ \ \ \ \ \ \ \ \ \ \ \ \ \ \ \ \ \ \ \ \ \ \ \ \ \ \ \ \ \ \ \ \ \ \ \ \ \ \ \ \ \ \ \ \ \ \ \ \ \ \ \ \ \ \ \ \ \ \ \ \ \ \ \ \ \ \ \ \ \ \ \ \ \ \ \ \ \ \ \ \ \ \ \ \ \ \ \ \
$$
\begin{equation}\label{0theorem1}
=e^{-\alpha n^{r}}n^{\frac{1-r}{s'}}\bigg(\frac{\|\cos t\|_{s}}{\pi^{1+\frac{1}{s}}(\alpha r)^{\frac{1}{s'}}} \mathcal{J}_{s}\Big(\frac{ \pi n^{1-r}}{\alpha r}\Big)+
\delta_{n,s}^{(1)}\Big(\frac{1}{(\alpha r)^{1+\frac{1}{s'}}} \mathcal{J}_{s}\Big(\frac{ \pi n^{1-r}}{\alpha r}\Big)\frac{1}{n^{r}}+\frac{1}{n^{\frac{1-r}{s'}}}\Big)\bigg),
\end{equation}
\begin{equation}\label{2theorem1}
\frac{1}{\pi}\inf\limits_{\lambda\in\mathbb{R}}\|P_{\alpha,r,\beta}^{(n)}(\cdot)-\lambda) \|_{s}= \ \ \ \ \ \ \ \ \ \ \ \ \ \ \ \ \ \ \ \ \ \ \ \ \ \ \ \ \ \ \ \ \ \ \ \ \ \ \ \ \ \ \ \ \ \ \ \ \ \ \ \ \ \ \ \ \ \ \ \ \ \ \ \ \ \ \ \ \ \ \ \ \ \ \ \ \
$$
$$
=e^{-\alpha n^{r}}n^{\frac{1-r}{s'}}\bigg(\frac{\|\cos t\|_{s}}{\pi^{1+\frac{1}{s}}(\alpha r)^{\frac{1}{s'}}} \mathcal{J}_{s}\Big(\frac{ \pi n^{1-r}}{\alpha r}\Big)+\delta_{n,s}^{(2)}\Big(\frac{1}{(\alpha r)^{1+\frac{1}{s'}}} \mathcal{J}_{s}\Big(\frac{ \pi n^{1-r}}{\alpha r}\Big)\frac{1}{n^{r}}+\frac{1}{n^{\frac{1-r}{s'}}}\Big)\bigg),
\end{equation}
$$
\frac{1}{2\pi}\sup\limits_{h\in\mathbb{R}}\|P_{\alpha,r,\beta}^{(n)}(\cdot)-P_{\alpha,r,\beta}^{(n)}(\cdot+h) \|_{s}= \ \ \ \ \ \ \ \ \ \ \ \ \ \ \ \  \ \ \ \ \ \ \ \ \ \ \ \ \ \ \ \ \ \ \ \ \ \ \ \ \ \ \ \ \ \ \ \ \ \ \ \ \ \ \ \ \ \ \ \ \ \
$$
\begin{equation}\label{1theorem1}
=e^{-\alpha n^{r}}n^{\frac{1-r}{s'}}\bigg(\frac{\|\cos t\|_{s}}{\pi^{1+\frac{1}{s}}(\alpha r)^{\frac{1}{s'}}} \mathcal{J}_{s}\Big(\frac{ \pi n^{1-r}}{\alpha r}\Big)+\delta_{n,s}^{(3)}\Big(\frac{1}{(\alpha r)^{1+\frac{1}{s'}}} \mathcal{J}_{s}\Big(\frac{ \pi n^{1-r}}{\alpha r}\Big)\frac{1}{n^{r}}+\frac{1}{n^{\frac{1-r}{s'}}}\Big)\bigg),
\end{equation}
де для величин ${\delta_{n,s}^{(i)}=\delta_{n,s}^{(i)}(\alpha,r,\beta)}$, $i=\overline{1,3}$ виконується нерівність  ${|\delta_{n,s}^{(i)}|\leq(14\pi)^{2}}$.

З формул (\ref{s12}), (\ref{0theorem1}) і (\ref{1theorem1}) випливає (\ref{theorem1}).
Теорему 1 доведено.

Наведемо деякі наслідки з теореми 1.

При $1<s<\infty$ з теореми 1 випливає наступне твердження.

{\bf Теорема 2.} {\it Нехай  $0<r<1$, $1<s<\infty$, $\alpha>0$ i $\beta\in\mathbb{R}$. Тоді  при
 $n\geq n_0(\alpha,r,s')$, $\frac{1}{s}+\frac{1}{s'}=1$, справедлива  оцінка
$$
 {\cal E}_{n}(C^{\alpha,r}_{\beta,1})_{s}=e^{-\alpha n^{r}}n^{\frac{1-r}{s'}}\bigg(
\frac{\|\cos t\|_{s}}{\pi^{1+\frac{1}{s}}(\alpha r)^{\frac{1}{s'}}}F^{\frac{1}{s}}\Big(\frac{1}{2}, \frac{3-s}{2}; \frac{3}{2}; 1\Big)+
$$
\begin{equation}\label{remark}
+
\gamma^{(2)}_{n,s}\Big(\Big(1+\frac{(\alpha r)^{\frac{s-1}{s'}}}{s-1}\Big) \frac{1}{n^{\frac{1-r}{s'}}}+\frac{(s')^{\frac{1}{s}}}{(\alpha r)^{1+\frac{1}{s'}}}\frac{1}{n^{r}}\Big)\bigg),
\end{equation}
де $F(a,b;c;z)$ --- гіпергеометрична функція Гаусса,
а  для величини ${\gamma_{n,s}^{(2)}=\gamma_{n,s}^{(2)}(\alpha,r,\beta)}$ виконується нерівність  ${|\gamma_{n,s}^{(2)}|\leq(14\pi)^{2}}$.}

{\bf Доведення теореми 2.}
Згідно з теоремою   1 для всіх $1<s<\infty$, ${0<r<1}$, $\alpha>0$, $\beta\in\mathbb{R}$ при $n\geq n_0(\alpha,r,s')$, $\frac{1}{s}+\frac{1}{s'}=1$, має місце наступна оцінка:
 $$
{\cal E}_{n}(C^{\alpha,r}_{\beta,1})_{s}=
e^{-\alpha n^{r}}n^{\frac{1-r}{s'}}\bigg(\frac{\|\cos t\|_{s}}{(\alpha r)^{\frac{1}{s'}}\pi^{1+\frac{1}{s}}}\bigg(\int\limits_{0}^{\frac{ \pi n^{1-r}}{\alpha r}}\frac{dt}{(t^2+1)^{\frac{s}{2}}}\bigg)^{\frac{1}{s}}+
$$
\begin{equation}\label{a4}
+\gamma_{n,s}^{(1)}\Big(\frac{1}{(\alpha r)^{1+\frac{1}{s'}}}\bigg(\int\limits_{0}^{\frac{ \pi n^{1-r}}{\alpha r}}\frac{dt}{(t^2+1)^{\frac{s}{2}}}\bigg)^{\frac{1}{s}}\frac{1}{n^{r}}+\frac{1}{n^{\frac{1-r}{s'}}}\Big)\bigg),
\end{equation}
де $\frac{1}{s}+\frac{1}{s'}=1$, і величина ${\gamma_{n,s}^{(1)}=\gamma_{n,s}^{(1)}(\alpha,r,\beta)}$ така, що  ${|\gamma_{n,s}^{(1)}|\leq(14\pi)^{2}}$.

 В [\ref{SerdyukStepanyuk2016}] доведено, що при $n\geq n_0(\alpha,r,s')$, виконується оцінка
 \begin{equation}\label{th2}
  \bigg(\int\limits_{0}^{\frac{ \pi n^{1-r}}{\alpha r}}\frac{dt}{(t^2+1)^{\frac{s}{2}}}\bigg)^{\frac{1}{s}}=\bigg(\int\limits_{0}^{\infty}\frac{dt}{(t^2+1)^{\frac{s}{2}}}\bigg)^{\frac{1}{s}}+
  \frac{\Theta^{(1)}_{\alpha,r,s',n}}{s-1}\Big(\frac{\alpha r}{\pi n^{1-r}}\Big)^{s-1},
 \end{equation}
 де  $|\Theta^{(1)}_{\alpha,r,s',n}|<2$.

 Легко бачити, що (див., наприклад, [\ref{grad}, с.962])
 \begin{equation}\label{da3}
\int\limits_{0}^{\infty}\frac{dt}{(t^2+1)^{\frac{s}{2}}}=\frac{1}{2}\int\limits_{0}^{\infty}\frac{t^{-\frac{1}{2}}dt}{(t+1)^{\frac{s}{2}}}=
 \frac{1}{2}B\Big(\frac{1}{2}, \frac{s-1}{2} \Big), \ \ \ 1<s<\infty,
 \end{equation}
 де
 $$
 B(x,y):=\int\limits_{0}^{\infty}\frac{u^{x-1}du}{(u+1)^{x+y}}, \ \ x, y>0,
 $$
 --- ейлеровий інтеграл 1-го роду.

Врахувавши  формулу [\ref{grad}, с. 964]
 $$
 B(x,y)=\frac{\Gamma(x)\Gamma(y)}{\Gamma(x+y)}=B(y,x),
 $$
 де
 $$
 \Gamma(z)=\int\limits_{0}^{\infty}t^{z-1}e^{-t}dt, \ \ \ Re z>0,
 $$
 --- ейлеровий інтеграл 2-го роду,
 та скориставшись теоремою Гаусса [\ref{grad}, с. 1056]
  $$
  F(a,b;c;1)=  \frac{\Gamma(c)\Gamma(c-a-b)}{\Gamma(c-a)\Gamma(c-b)}
 $$
 при  $a=\frac{1}{2}$, $b=\frac{3-s}{2}$, $c=\frac{3}{2}$, з (\ref{da3}) отримуємо
  $$
\int\limits_{0}^{\infty}\frac{dt}{(t^2+1)^{\frac{s}{2}}}= \frac{1}{2}\frac{\Gamma\big( \frac{1}{2}\big)\Gamma\big( \frac{s-1}{2}\big)}{\Gamma\big( \frac{s}{2}\big)}
   =\frac{\sqrt{\pi}}{2}\frac{\Gamma\big( \frac{s-1}{2}\big)}{\Gamma\big( \frac{s}{2}\big)}=\frac{\Gamma\big( \frac{3}{2}\big)\Gamma\big( \frac{s-1}{2}\big)}{\Gamma\big( \frac{s}{2}\big)} =F\Big(\frac{1}{2}, \frac{3-s}{2}; \frac{3}{2}; 1\Big).
 $$
 Отже, для довільних $1<s<\infty$ має місце рівність
  \begin{equation}\label{d1}
\bigg(\int\limits_{0}^{\infty}\frac{dt}{(t^2+1)^{\frac{s}{2}}}\bigg)^{\frac{1}{s}}= F^{\frac{1}{s}}\Big(\frac{1}{2}, \frac{3-s}{2}; \frac{3}{2}; 1\Big).
 \end{equation}

 Врахувавши оцінку
 \begin{equation}\label{th21}
\bigg(\int\limits_{0}^{\frac{ \pi n^{1-r}}{\alpha r}}\frac{dt}{(t^2+1)^{\frac{s}{2}}}\bigg)^{\frac{1}{s}}\leq
\bigg(\int\limits_{0}^{\infty}\frac{dt}{(t^2+1)^{\frac{s}{2}}}\bigg)^{\frac{1}{s}} <
 \bigg(1+\int\limits_{1}^{\infty}\frac{dt}{t^{s}}\bigg)^{\frac{1}{s}}<
(s')^{\frac{1}{s}},
 \end{equation}
з формул (\ref{a4}), (\ref{th2}), (\ref{d1}) i (\ref{th21}) отримуємо
$$
 {\cal E}_{n}(C^{\alpha,r}_{\beta,1})_{s}=e^{-\alpha n^{r}}n^{\frac{1-r}{s'}}\bigg(
\frac{\|\cos t\|_{s}}{\pi^{1+\frac{1}{s}}(\alpha r)^{\frac{1}{s'}}}F^{\frac{1}{s}}\Big(\frac{1}{2}, \frac{3-s}{2}; \frac{3}{2}; 1\Big)+
$$
\begin{equation}\label{remark0}
+
\gamma^{(3)}_{n,s}\Big(\frac{1}{s-1}\frac{(\alpha r)^{\frac{s-1}{s'}}}{n^{(1-r)(s-1)}}+\frac{(s')^{\frac{1}{s}}}{(\alpha r)^{1+\frac{1}{s'}}}\frac{1}{n^{r}}+\frac{1}{n^{\frac{1-r}{s'}}}\Big)\bigg),
\end{equation}
де   для величини ${\gamma_{n,s}^{(3)}=\gamma_{n,s}^{(3)}(\alpha,r,\beta)}$ виконується нерівність  ${|\gamma_{n,s}^{(3)}|\leq(14\pi)^{2}}$.

Оскільки
$$
\frac{1}{n^{\frac{1-r}{s'}}}=\frac{1}{n^{(1-r)\frac{s-1}{s}}}>\frac{1}{n^{(1-r)(s-1)}},
$$
то
\begin{equation}\label{remarkk}
\frac{1}{s-1}\frac{(\alpha r)^{\frac{s-1}{s'}}}{n^{(1-r)(s-1)}}+\frac{1}{n^{\frac{1-r}{s'}}}
\leq\Big(1+\frac{(\alpha r)^{\frac{s-1}{s'}}}{s-1}\Big) \frac{1}{n^{\frac{1-r}{s'}}}.
\end{equation}
З (\ref{remark0}) і (\ref{remarkk}) випливає (\ref{remark}).
Теорему 2 доведено.

З теореми 2 при  $s=2$ одержуємо наступне твердження.

{\bf Наслідок 1.} {\it Нехай  $0< r<1$,  $\alpha>0$ і  $\beta\in\mathbb{R}$. Тоді при  $n\geq n_0(\alpha,r,2)$ справедлива  оцінка
\begin{equation}\label{consequence2}
{\cal E}_{n}(C^{\alpha,r}_{\beta,1})_{2}
 =\frac{e^{-\alpha n^{r}}}{\sqrt{2\pi\alpha r}}n^{\frac{1-r}{2}}\bigg(1+
\gamma^{(1)}_{n}\Big(\Big(\alpha r+\sqrt{\alpha r}\Big)\frac{1}{n^{\frac{1-r}{2}}}
 +\frac{\sqrt{2}}{\alpha r}\frac{1}{n^{r}}\Big)\bigg).
\end{equation}
де  для величини ${\gamma_{n}^{(1)}=\gamma_{n}^{(1)}(\alpha,r,\beta)}$ виконується нерівність  ${|\gamma_{n}^{(1)}|\leq392\pi^{\frac{5}{2}}}$.}

{\bf Доведення наслідку 1.} Дійсно, поклавши в рівності (\ref{remark}) $s=s'=2$, врахувавши, що $F\Big(\frac{1}{2},\frac{1}{2};\frac{3}{2};1\Big)=\frac{\pi}{2}$, при  $n\geq n_0(\alpha,r,2)$  отримуємо
$$
 {\cal E}_{n}(C^{\alpha,r}_{\beta,1})_{2}=e^{-\alpha n^{r}}n^{\frac{1-r}{2}}\bigg(
\frac{\|\cos t\|_{2}}{\pi^{\frac{3}{2}}(\alpha r)^{\frac{1}{2}}}\bigg(F^{\frac{1}{2}}\Big(\frac{1}{2}, \frac{1}{2}; \frac{3}{2}; 1\Big)+
\gamma^{(2)}_{n,2}\Big(\Big(1+\sqrt{\alpha r}\Big) \frac{1}{n^{\frac{1-r}{2}}}
 +\frac{\sqrt{2}}{(\alpha r)^{\frac{3}{2}}}\frac{1}{n^{r}}\Big)\bigg)=
$$
$$
=\frac{e^{-\alpha n^{r}}}{\sqrt{2\pi\alpha r}}n^{\frac{1-r}{2}}\bigg(
1+
\gamma^{(2)}_{n,2}\sqrt{2\pi}\Big(\Big(1+\sqrt{\alpha r}\Big) \frac{\sqrt{\alpha r}}{n^{\frac{1-r}{2}}}
 +\frac{\sqrt{2}}{\alpha r}\frac{1}{n^{r}}\Big)\bigg).
$$
Наслідок 1 доведено.

Утім більш точну, ніж (\ref{consequence2}), оцінку можна одержати  виходячи
з рівності (\ref{serd2011}). А саме, з урахуванням (\ref{serd2011}), а також   формули (10) з роботи [\ref{SerdyukStepanyuk2016dop}] при $\alpha>0$, $r\in(0,1)$, $\beta\in\mathbb{R}$ i $n\geq n_0(\alpha,r,2)$ має місце оцінка
\begin{equation}\label{sta1}
{\cal E}_{n}(C^{\alpha,r}_{\beta,1})_{2}=\frac{1}{\sqrt{\pi}}\Big(\sum\limits_{k=n}^{\infty}e^{-2\alpha k^{r}}\Big)^{\frac{1}{2}}=\frac{e^{-\alpha n^{r}}}{\sqrt{2\pi\alpha r}}n^{\frac{1-r}{2}}\Big(
1+\gamma_{n}^{(2)}\Big(\frac{1}{2\alpha r}\frac{1}{n^{r}}+\frac{\alpha r}{n^{1-r}}\Big)\Big),
\end{equation}
в якій  для величини ${\gamma_{n}^{(2)}=\gamma_{n}^{(2)}(\alpha,r)}$ виконується нерівність  $|\gamma_{n}^{(2)}|\leq\sqrt{\frac{54\pi^{3}}{54\pi^{3}-1}} $.

При $s=1$ теорема 1 дозволяє  уточнити асимптотичну рівність (\ref{step}).
Має місце наступне твердження.

{\bf Теорема 3.} {\it Нехай  $0< r<1$,  $\alpha>0$  i  $\beta\in\mathbb{R}$. Тоді при  $n\geq n_1(\alpha,r)$ справедлива  оцінка
\begin{equation}\label{consequence}
{\cal E}_{n}(C^{\alpha,r}_{\beta,1})_{1}=\frac{4}{\pi^{2}}e^{-\alpha n^{r}}\ln\Big(\frac{\pi n^{1-r}}{\alpha r}\Big)+\gamma_{n,1}^{(2)}e^{-\alpha n^{r}},
\end{equation}
де  для величини ${\gamma_{n,1}^{(2)}=\gamma_{n,1}^{(2)}(\alpha,r,\beta)}$ виконується нерівність  ${|\gamma_{n,1}^{(2)}|\leq20\pi^{4}}$.}

{\bf Доведення теореми 3.} Із (\ref{n_p}) i (\ref{n_2}) випливає, що ${n_1(\alpha,r)>n_0(\alpha,r,1)}$. Тому, поклавши в рівності (\ref{theorem1}) $s=1$ та врахувавши оцінку
$$
\int\limits_{0}^{\frac{\pi n^{1-r}}{\alpha r}}\frac{dt}{\sqrt{t^{2}+1}}=
\int\limits_{1}^{\frac{\pi n^{1-r}}{\alpha r}}\frac{dt}{t}+
\bigg(\int\limits_{0}^{\frac{\pi n^{1-r}}{\alpha r}}\frac{dt}{\sqrt{t^{2}+1}}-\int\limits_{1}^{\frac{\pi n^{1-r}}{\alpha r}}\frac{dt}{t}\bigg) =\ln\Big(\frac{\pi n^{1-r}}{\alpha r}\Big)+\Theta_{\alpha,r,n}, \ \ \ 0<\Theta_{\alpha,r,n}<1,
$$
при  $n\geq n_1(\alpha,r)$ отримаємо
$$
{\cal E}_{n}(C^{\alpha,r}_{\beta,1})_{1} =
  e^{-\alpha n^{r}}\bigg(
\frac{4}{\pi^{2}}\int\limits_{0}^{\frac{\pi n^{1-r}}{\alpha r}}\frac{dt}{\sqrt{t^{2}+1}}
+\gamma_{n,\infty}^{(1)}\bigg(\frac{1}{\alpha r}\frac{1}{n^{r}}\int\limits_{0}^{\frac{\pi n^{1-r}}{\alpha r}}\frac{dt}{\sqrt{t^{2}+1}}+1\bigg)\bigg)=
$$
\begin{equation}\label{1f111}
=
  e^{-\alpha n^{r}}\bigg(
\frac{4}{\pi^{2}}\ln\Big(\frac{\pi n^{1-r}}{\alpha r}\Big)+\frac{4}{\pi^{2}}\Theta_{\alpha,r,n}
+\gamma_{n,\infty}^{(1)}\bigg(\frac{1}{\alpha r}\frac{1}{n^{r}}\ln\Big(\frac{\pi n^{1-r}}{\alpha r}\Big)+\frac{\Theta_{\alpha,r,n}}{\alpha rn^{r}}+1\bigg)\bigg).
 \end{equation}
Підрахунки показують, що  при  $n\geq n_1(\alpha,r)$
\begin{equation}\label{2f111}
 \frac{4}{\pi^{2}}\Theta_{\alpha,r,n}+
|\gamma_{n,\infty}^{(1)}|\bigg(\frac{1}{\alpha r}\frac{1}{n^{r}}\ln\Big(\frac{\pi n^{1-r}}{\alpha r}\Big)+\frac{\Theta_{\alpha,r,n}}{\alpha rn^{r}}+1\bigg)\leq 20\pi^{4},
 \end{equation}
а, тому з (\ref{1f111}) i (\ref{2f111}) отримуємо (\ref{consequence}).
Теорему 3 доведено.

Зі співвідношення (\ref{consequence}) випливає встановлена  О.І. Степанцем асимптотична рівність (\ref{step}).

На завершення роботи зауважимо, що з урахуванням формул (\ref{d1}) i (\ref{remarkk}) один з основних результатів роботи [\ref{SerdyukStepanyuk2016dop}] --- оцінка (8) з теореми 2 --- може бути сформульоваим наступним чином.

{\bf Теорема 4.} {\it Нехай  $0<r<1$, $1<p<\infty$, $\alpha>0$ i $\beta\in\mathbb{R}$. Тоді  при
 $n\geq n_0(\alpha,r,p)$ справедлива  оцінка
$$
 {\cal E}_{n}(C^{\alpha,r}_{\beta,p})_{\infty}=e^{-\alpha n^{r}}n^{\frac{1-r}{p}}\bigg(
\frac{\|\cos t\|_{p'}}{\pi^{1+\frac{1}{p'}}(\alpha r)^{\frac{1}{p}}}F^{\frac{1}{p'}}\Big(\frac{1}{2}, \frac{3-p'}{2}; \frac{3}{2}; 1\Big)+
$$
\begin{equation}\label{remark}
+
\gamma^{(2)}_{n,p}\Big(\Big(1+\frac{(\alpha r)^{\frac{p'-1}{p}}}{p'-1}\Big) \frac{1}{n^{\frac{1-r}{p}}}+\frac{(p)^{\frac{1}{p'}}}{(\alpha r)^{1+\frac{1}{p}}}\frac{1}{n^{r}}\Big)\bigg),
\end{equation}
в якій $\frac{1}{p}+\frac{1}{p'}=1$, $F(a,b;c;z)$ --- гіпергеометрична функція Гаусса,
а  для величини ${\gamma_{n,p}^{(2)}=\gamma_{n,p}^{(2)}(\alpha,r,\beta)}$ виконується нерівність  ${|\gamma_{n,p}^{(2)}|\leq(14\pi)^{2}}$.}
\vskip 10mm

\begin{enumerate}

\item \label{SerdyukStepanyuk2016dop}
{\it Сердюк А.С., Степанюк Т.А. }
Рівномірні наближення сумами Фур'є на класах згорток з узагальненими ядрами Пуассона // Доповіді НАНУ --- 2016. --- №11. --- С. 10--16.

\item \label{SerdyukStepanyuk2016}
{\it Serdyuk A.S., Stepaniuk T.A., }
Uniform approximations by Fourier sums on classes of generalized Poisson integrals // Arxiv preprint,  arXiv:1603.01891, (2016), 31 p.

\item \label{Stepanets1}
{\it Степанец А.И.} Методы теории приближений: В 2 ч. // Праці Інституту математики НАН України --- Киев: Ин-т
математики НАН Украины, 2002. --- {\bf 40}. --- Ч.І. --- 427 с.

\item\label{Stepanets_Serdyuk_Shydlich} {\it Степанець О.І., Сердюк А.С., Шидліч А.Л.}  Про деякі нові критерії нескінченної диференційовності періодичних функцій
// Укр. мат. журн. --- 2007.
--- {\bf 59}, №10. --- С.~1399--1409.

\item\label{Stepanets_Serdyuk_Shydlich2009} {\it Степанец А.И., Сердюк А.С., Шидлич А.Л.} О связи классов $(\psi,\beta)$--дифференцируемых функций с классами Жевре // Укр. мат. журн. --- 2009. --- {\bf 61}, № 1. --- С. 140-145.

\item\label{Stepanets2}{Степанец А.И.\/} Методы теории
приближений: В 2 ч.   // Пр. Iн-ту математики НАН України T. 40. --- К.\,: Ин-т
математики НАН Украины, 2002. --- Ч. ІI.
--- 427 с.

\item\label{S_S}
{\it Сердюк А.С., Степанюк Т.А.}  Порядкові оцінки найкращих
наближень і наближень сумами Фур'є класів нескінченно
диференційовних функцій// Теорія наближення функцій та суміжні питання: Зб. праць Ін--ту математики НАН України. --- Київ: Інститут математики НАН України,
2013. --- {\bf 10}, №1.
 --- С. 255--282.

\item\label{S_S2}
{\it Сердюк А.С., Степанюк Т.А.}  Оцінки найкращих наближень класів нескінченно диференційовних функцій
в рівномірній та інтегральній метриках  // Укр. мат. журн. --- 2014. ---  {\bf 66}, №9. ---  С.1244--1256.

\item \label{SerdyukStepanyukBulleten}
{\it A.S. Serdyuk and T.A. Stepanyuk,} Estimates for approximations by Fourier sums, best approximations and best orthogonal trigonometric approximations of the classes of $(\psi,\beta)$--differentiable functions, Bulletin: de la societe des sciences et des letters de Lodz, vol. \textbf{66}, No 2 (2016).

    \item \label{Nikolsky 1946}
{\it Никольский С.М.} Приближение функций тригонометрическими полиномами в среднем
// Известия АН СССР. Сер. мат. --- 1946. --- {\bf 10}, №3.
 --- С. 207--256.

\item \label{Stechkin 1980}
{\it Стечкин С.Б.} Оценка остатка ряда Фурье для дифференцируемых функций
// Труды Мат. Ин-та АН СССР. --- 1980. --- {\bf 145}.
 --- С. 126--151.

\item \label{Serdyuk2005Int}
{\it Сердюк А.С.} Наближення класів аналітичних функцій сумами Фур’є в метриці простору $L_{p}$
// Укр. мат. журн. --- 2005. --- {\bf 57}, №10.
 --- С. 1395–-1408.

\item \label{Serdyuk2012}
{\it
Сердюк А.С.} Наближення інтерполяційними тригонометричними поліномами на класах періодичних аналітичних функцій // Укр. мат. журн. --- 2012. --- {\bf 64}, № 5. --- С.~698--712.

\item \label{Serdyuk2013}
{Сердюк А.С., Соколенко І.В.\/} Наближення лінійними методами
класів  $(\psi, \overline{\beta})$--диференційовних функцій
//  Теорія наближення функцій та суміжні питання: Зб. праць Ін--ту математики НАН України. ---  2013. --- Т.~10, №1.
 --- С. 245--254.

\item \label{Step monog 1987} {\it  Степанец А.И.} Классификация и приближение периодических функций. ---
Киев: Наук. думка~--- 1987.~--- 268 c.

\item \label{Teljakovsky 1989}
{\it Теляковский С.А.} О приближении суммами Фурье функций высокой гладкости // Укр. матем. журн.--- 1989. --- {\bf 41}, №4.
 --- С. 510--518.

\item\label{grad}
Грандштейн И.С., Рыжик И.М. Таблицы интегралов,
сумм, рядов и произведений. --- М.: Гос. из-во физ.-мат. лит-ры, 1963. --- 1100 с.

\end{enumerate}

\end{document}